\documentclass[12pt,doublespace]{article}

\newcommand{\pr}{{\it Proof.} \ }
\newcommand{\e}{\alpha}
\newcommand{\f}{\beta_{t}}

\usepackage{amsfonts}
\usepackage{amsmath}
\usepackage{amssymb}
\usepackage{latexsym}

\begin{document}
\title{The fourth moment of $\zeta^{'}(\rho)$}
\author{Nathan Ng}
\maketitle

\begin{abstract}
Discrete moments of the Riemann zeta function were studied
by Gonek and Hejhal in the 1980's.  They independently formulated
a conjecture concerning the size of these moments. In 1999, Hughes,
Keating, and O'Connell, by employing a random matrix model, made this
conjecture more precise.
Subject to the Riemann hypothesis,
we establish upper and lower bounds of
the correct order of magnitude in the case of the fourth moment.
\end{abstract}

\section{Introduction}

\footnotetext{2000 Mathematics Subject Classification. Primary 11026;
 Secondary 11M26.}

This article concerns discrete moments of the derivative of the
Riemann zeta function of the form
\[ J_{k}(T) = \sum_{0 < \gamma \le T} |\zeta^{'}(\rho)|^{2k}
\]
where $\rho = \beta + i \gamma$ ranges over non-trivial zeros
of $\zeta(s)$ and $k \in \mathbb{R}$.
In particular, we focus on the case $k=2$.
These moments are discrete analogues of the ordinary moments
of the Riemann zeta function.
In recent years there has been renewed interest in the moments of
$L$-functions, in part due to Keating and Snaith's \cite{KS}
work in random matrix theory.
Estimates for the discrete moments have number theoretic applications
(see  \cite{CGG},\cite{Mu},\cite{Ng}).
To date, few asymptotic formulae have been
established
for these moments. However,
Gonek \cite{G3} and Hejhal \cite{He}
independently conjectured
\begin{equation}
   J_{k}(T) \asymp T \log^{(k+1)^{2}} T
\end{equation}
for $k \in \mathbb{R}$.
Hughes, Keating,
and O'Connell \cite{HKO},
applying random matrix models refined this to: \\\\
{\bf Random Matrix Model Conjecture}
For $k > -\frac{3}{2}$ and bounded, \\
\begin{equation}
   J_{k}(T) \sim \frac{G^{2}(k+2)}{G(2k+3)} \cdot a_{k} \cdot N(T)
                 \cdot \left( \log \frac{T}{2\pi} \right)^{k(k+2)}
   \label{eq:rmmc}
\end{equation}
as $T \rightarrow \infty$,where $G$ is Barnes' function defined by
\[ G(z+1) =(2 \pi)^{z/2}\exp \left(-\frac{1}{2}(z^{2} + \gamma z^{2} + z)
   \right)
   \prod_{n=1}^{\infty}
   \left(
   \left(1 + \frac{z}{n} \right)^{n} e^{-z + z^{2}/2n}
   \right) \ ,
\]
$\gamma$ is Euler's constant,
$  a_{k} = \prod_{p} \left( 1 - \frac{1}{p} \right)^{k^{2}}
          \sum_{m = 0}^{\infty} \left( \frac{\Gamma(m+k)}{m! \Gamma(k)}
          \right)^{2}p^{-m}$,
and $N(t)$ denotes the number of zeros of $\zeta(s)$ in the box
with vertices $0,1,1+it,it$. \\

The number  $a_{2} = \zeta(2)^{-1} = \frac{6}{\pi^2}$
appears frequently in this article.
Conjecture~(\ref{eq:rmmc}) agrees with
results of Von Mangoldt and Gonek \cite{G1} in the cases $k = 0,1$.
Furthermore, one verifies
$J_{-1}(T) \sim \frac{3}{\pi^{3}} T$ is the case $k=-1$.
Gonek first conjectured
this formula by methods similar to Montgomery's study of the
pair correlation conjecture.
When $k=2$,~(\ref{eq:rmmc}) reduces to
$J_{2}(T) \sim  \frac{1}{2880 \pi^{3}} T \log^{9} T$.
We establish that the random matrix theory conjecture
is of the correct order of magnitude in this case.
Throughout, we use the notation
$L = \log \frac{T}{2 \pi}$.
Our main result is
\newtheorem{4}{Theorem}
\begin{4}
The Riemann hypothesis implies
\begin{equation}
  \frac{c_{1}}{\pi^{3}}
 TL^{9} \left(1+O\left(\frac{\log L}{L} \right) \right) \le
  J_{2}(T)
  \le \frac{c_{2}}{\pi^{3}}
  TL^{9} \left(1+O\left(\frac{\log L}{L}\right) \right)
\end{equation}
\end{4}
where
\begin{equation}
  c_{1}= (\sqrt{a}-\sqrt{b})^{2} = 0.0000687... \ , \
  c_{2}=(\sqrt{a}+\sqrt{b})^{2} = 0.0051561...
  \label{eq:constants}
\end{equation}
with $a = \frac{61}{60480},b=\frac{97}{60480}$.
In contrast, $\frac{1}{2880} = 0.0003472...$. \\

The same techniques as Theorem 1, permit one to replace $\zeta^{'}(s)$
by higher derivatives.  We remark that
only Theorem 1 depends on RH.  All other lemmas, corollaries,
and theorems are independent of any hypothesis.
We establish the following unconditional result which
may be of use in future moment calculations.
\newtheorem{smoo}[4]{Theorem}
\begin{smoo}
Let $d(n)$ denote the number of divisors of $n$ and $\delta =
\lambda/ \log \left( \frac{T}{2 \pi} \right)$ where
$\lambda \in \mathbb{R}$ and $|\lambda| \ll 1$. Then we have
\begin{equation}
\begin{split}
   & \sum_{0 < \gamma < T} \sum_{m \le \frac{\gamma}{2 \pi}}
   \frac{d(m)}{m^{\rho+i\delta}}
   \sum_{n \le \frac{\gamma}{2 \pi}}
   \frac{d(n)}{n^{1-\rho- i\delta}}
    = \frac{3}{\pi^{3}}
    \left( \frac{1}{5!} - 4 \sum_{j \ge 1} \frac{(-1)^{j}
    \lambda^{2j}}{(5+2j)!} \right)
    TL^{5}(1+o(1))   \\
\end{split}
\end{equation}
where $\rho = \beta + i \gamma$ ranges over non-trivial zeros of
the zeta function with $0 < \gamma <T$. The $o(1)$ term is $(\log
L)/L$.
\end{smoo}
{\bf Notation}  We work with Dirichlet series of the form
\begin{equation}
   (-1)^{\mu+\nu} \zeta^{(\mu)}(s) \zeta^{(\nu)}(s) =
   \sum_{n=1}^{\infty} \frac{d^{(\mu,\nu)}(n)}{n^{s}}
   \label{eq:duvc}
\end{equation}
where $\mu, \nu \in \mathbb{Z}_{\ge 0}$. Note that $d^{(\mu,
\nu)}(n) = (\log^{\mu}* \log^{\nu})(n)$ where $*$ denotes
convolution. Furthermore, we set $d^{(\mu)}(n) := d^{(\mu,0)}(n)$.
The generalized divisor function $d_{k}(n)$ for $k > 0$ is defined
by its generating function $\zeta^{k}(s)= \sum_{n=1}^{\infty}
\frac{d_{k}(n)}{n^{s}}$. In this article the arithmetic functions
\begin{equation}
\begin{split}
    \e(n) & := d^{(1,1)}(n)
     = (\log n) \, d^{(1)}(n) - d^{(2)}(n)  \ , \\
    \f(n) & :=  (a_{t}*a_{t})(n)
     = l^{2} d(n)
   - 2 l \ d^{(1)}(n) + \e(n) \ ,
    \label{eq:froid} \\
\end{split}
\end{equation}
where $l = \log( \frac{t}{2 \pi})$
and $a_{t}(n) = \log( \frac{t}{2 \pi n})$
appear often.
To simplify notation, we define for an
arbitrary sequence $a(n,t)$ with $n \in \mathbb{Z}^{+}$
and $t \in \mathbb{R}$  the Dirichlet polynomial
\begin{equation}
   D_{a}(\sigma + it) = \sum_{n \le \frac{t}{2\pi}}
   \frac{a(n,t)}{n^{\sigma +i t}}
   \ .
   \label{eq:dirp}
\end{equation}
{\bf Acknowledgements}  The author thanks Professor Andrew
Granville for helpful discussions concerning this article.

\subsection{Proof of Theorem 1}
\indent We commence with the proof of Theorem 1 since the the
argument is rather simple. This proof is subject to Corollary 1, a
mean value result,  which is a special case of Lemma 5. However,
the proofs of Corollary 1 and Theorem 2 are deferred until later.
We first state Corollary 1.
\newtheorem{mvt}{Corollary}
\begin{mvt}
We have
\begin{equation}
    S_{\alpha} = \sum_{0 < \gamma \le T}
     D_{\e}(\rho) D_{\e}(1-\rho)
      =  \frac{61}{60480 \pi^{3}}  T L^9
     + O(T L^{8} \log L) \ ,
     \label{eq:sa}
\end{equation}
\begin{equation}
    S_{\beta} = \sum_{0 < \gamma \le T}
    D_{\beta_{\gamma}}(\rho) D_{\beta_{\gamma}}(1-\rho)
    = \frac{97}{60480 \pi^{3}}
     T L^{9} + O(T L^{8} \log L)
    \label{eq:sb}
\end{equation}
where $\rho = \beta + i \gamma$ ranges through the non-trivial zeros
of the zeta function with $0 < \gamma < T$.
Note that $D_{\alpha}(s)$ and $D_{\beta_{\gamma}}(s)$ are Dirichlet
polynomials associated to $\alpha(n)$ and $\beta_{\gamma}(n)$
as defined by~(\ref{eq:dirp}).
\end{mvt}
{\it Proof of Theorem 1.} \,
The approximate functional equation we require is
\begin{equation}
  \zeta^{'}(\sigma + it)^{2} = \sum_{n \le \frac{|t|}{2 \pi}}
   \frac{\e(n)}{n^{\sigma + it}}  + \chi^{2}(\sigma + it)
   \sum_{n \le \frac{|t|}{2 \pi}}
   \frac{\beta_{t}(n)}{n^{1-\sigma - it}} + O( \log^{3} t )
   \label{eq:afe}
\end{equation}
where $\alpha(n)$ and $\beta_{t}(n)$ are defined by~(\ref{eq:froid}) and $\chi(s) = \pi^{s-\frac{1}{2}}
\Gamma(\frac{1-s}{2}) / \Gamma(\frac{s}{2}) $
is the factor from the functional equation of the zeta function.
It satisfies $\zeta(s) = \chi(s) \zeta(1-s) \ \mathrm{and} \
   \chi(s) \chi(1-s) = 1$.
Equation~(\ref{eq:afe}) is derived in \cite{C} (Lemma 3 p.29).
Let $\rho$ denote a non-trivial zero of the Riemann zeta
function.  By~(\ref{eq:afe}) we have
\begin{equation}
\begin{split}
  \zeta^{'}(\rho)^{2}
  \zeta^{'}(1-\rho)^{2} =
  & (D_{\alpha}(\rho) + \chi^{2}(\rho)
      D_{\beta_{\gamma}}(1-\rho) + O(l^{3}))
     \cdot \\
    &   (D_{\alpha}(1-\rho) + \chi^{2}(1-\rho)
      D_{\beta_{\gamma}}(\rho) + O(l^{3})) \\
  \label{eq:mieux}
\end{split}
\end{equation}
where $l = \log \gamma$. Summing~(\ref{eq:mieux}) over zeros that satisfy
$0 < \mathrm{Im}(\rho) < T$ yields
\begin{equation}
  \sum_{0 < \gamma < T}
  \zeta^{'}(\rho)^{2}
  \zeta^{'}(1-\rho)^{2}
  = S_{1} + 2 \mathrm{Re}(S_{2}) + S_{3} + S_{4}
 \label{eq:ide}
\end{equation}
where $S_{1} = S_{\alpha}+S_{\beta}$, $
   S_{2} =
    \sum_{0 < \gamma < T}
    \chi^{2}(1-\rho) D_{\e}(\rho)
    D_{\beta_{\gamma}}(\rho)$,
\[
    S_{3} \ll L^{3}
    \sum_{0 < \gamma < T}
    \left( \left| D_{\e}(\rho) \right| +
    \left|\chi^{2}(1-\rho)D_{\beta_{\gamma}}(\rho) \right|
    \right) \ ,
\]
and $S_{4} \ll (\log^{6} T) N(T) \ll T L^{7}$. We have
by Corollary 1
\[
  S_{1} = S_{\alpha}+S_{\beta}
        = \frac{a+b}{\pi^3}T L^{9}
   + O(T L^{8} \log L)
\]
where $a = \frac{61}{60480}$ and $b = \frac{97}{60480}$.
Note that under the assumption of RH $|\chi(1-\rho)|=1$ and
$|D_{a}(\rho)|^{2} = D_{a}(\rho)D_{a}(1-\rho)$  for a real
sequence $a=a(n,t)$. Hence assuming RH, Cauchy-Schwarz implies
\[
  |S_{2}|  \le S_{\alpha}^{\frac{1}{2}} S_{\beta}^{\frac{1}{2}}
          = \frac{\sqrt{a b}}{\pi^3}  T L^{9}
              \left( 1 + O(L^{-1} \log L) \right)
\]
and also
$S_{3} \ll  (N(T)L^{6})^{\frac{1}{2}}S_{1}^{\frac{1}{2}} \ll
   T L^{\frac{7}{2}} S_{1}^{\frac{1}{2}}$.  Lastly we note that
RH implies
\begin{equation}
  |\zeta^{'}(\rho)|^{4} = \zeta^{'}(\rho)^{2}
  \zeta^{'}(1-\rho)^{2} \ .
  \label{eq:absz}
\end{equation}
By~(\ref{eq:ide}),~(\ref{eq:absz}), and collecting our estimates
of the $S_{i}$ for $i=1,\ldots, 4$ we have
\[
   \frac{c_{1}}{\pi^3} T L^9(1 + O(L^{-1} \log L))
   \le J_{2}(T)
   \le
   \frac{c_{2}}{\pi^3}
    T L^9(1 + O( L^{-1} \log L)) \\
\]
for $c_{1},c_{2}$ as in~(\ref{eq:constants})
and Theorem 1 is established. \\\\
\indent In the above calculation RH was used to evaluate $S_{2}$
and $S_{3}$ and to guarantee the identity~(\ref{eq:absz}). It may
be possible, by more sophisticated techniques, to bound $S_{2}$
and $S_{3}$ independent of RH and obtain unconditional bounds for
the sum in~(\ref{eq:ide}). Moreover, we expect $S_{2}$ to
contribute to the main term of  $J_{2}(T)$. In contrast, the
analogous sum in Ingham's \cite{In} calculation does not
contribute.
\section{Lemmas}
\indent Our calculations require an old formula
of Landau's. We apply Gonek's uniform version
(proven in \cite{G2} pp.401-403).
\newtheorem{gonek}{Lemma}
\begin{gonek} \ Let $x,T > 1$ then \\
\begin{equation}
\begin{split}
  \sum_{0 < \gamma \le T} x^{\rho}
  & = - \frac{T}{2 \pi} \Lambda(x) +
    O \left( x (\log(2xT)) ( \log \log 3x) \right) \\
  & O \left( (\log x)
    \min \left( T, \frac{x}{\langle x \rangle} \right)
    \right) +
    O \left( (\log (2T)) \min \left( T, \frac{1}{\log x} \right)
    \right) \\
    \label{eq:gl}
\end{split}
\end{equation}
where $\langle x \rangle$
denotes the distance from $x$ to the nearest
prime power other than $x$ itself.
\end{gonek}

To prove Lemmas 3 and 5
we require estimates for divisor sums.
We only need upper bounds for shifted divisor sums as in $(i)$ below.
Moreover, we do not require
the stronger asymptotic formulae that have been proven.  In addition,
a Brun-Titchmarsh result for divisor sums is applied.
\newtheorem{ing}[gonek]{Lemma}
\begin{ing}
$(i)$ If $r \le x$ is a positive integer and $\sigma_{-1}(r) =
\sum_{d \mid r} d^{-1}$ then
\begin{equation}
   \sum_{r < n \le x} d(n) d(n-r) \ll \sigma_{-1}(r) x \log^2 x
   \ .
   \label{eq:shift}
\end{equation}
$(ii)$ Let $\lambda \in \Bbb R$, $k \in \Bbb N$, $a \in \Bbb Z$,
$(a,k) = 1$ and $k < x^{1-\alpha}$ for any $\alpha
>0$, then
\begin{equation}
   \sum_{{\begin{substack}{n \le x
         \\ n \equiv a \ \mathrm{mod} \ k}\end{substack}}}
   d_{r}^{\lambda}(n) \ll
   \frac{x}{k} \left( \frac{\phi(k)}{k} \log x
   \right)^{r^{\lambda}-1} \ .
   \label{eq:vino}
\end{equation}
\end{ing}
\pr Part $(i)$ is Lemma B2 of \cite{In} p.296 and
part $(ii)$ is a direct application of
Theorem 2 of \cite{Sh} p.169. \\

We prove a general mean value result for sequences which behave
like $d(n)$. Extending the following result to $d_{k}(n)$ for $k
\ge 3$ would require knowledge of sums like~(\ref{eq:shift}) with
$d(n)$ replaced by $d_{k}(n)$. However, such results have not been
proven yet.
\newtheorem{s1}[gonek]{Lemma}
\begin{s1}
Suppose two sequences $a(n)$and $b(n)$ satisfy $a(n) \ll
\log^{A}(n) d(n)$ and $b(n) \ll \log^{B}(n) d(n)$ for $A,B >0$.
Then we define for $\delta \in \mathbb{R}$ the mean values
\begin{equation}
\begin{split}
  I=I(a,b;T,\delta) & =
  \sum_{0 < \gamma < T} D_{a} \left( \rho+i \delta \right)
                        D_{b} \left( 1- \rho - i \delta
  \right)  \ , \\
  I(a,b;T) & := I(a,b;T,0)  \\
  \label{eq:iabt}
\end{split}
\end{equation}
and we have
\begin{equation}
\begin{split}
  I(a,b;T,\delta)  = \frac{T}{2 \pi} & \left(
   \log \left(\frac{T}{2 \pi} \right)
   \sum_{n \le \frac{T}{2 \pi}} \frac{a(n)b(n)}{n} -
      \sum_{mj \le \frac{T}{2 \pi}} \frac{\Lambda(j)a(m)b(mj)}{mj^{1-i\delta}}
      \right. \\
   &  \left. -  \sum_{mj \le \frac{T}{2 \pi}}
   \frac{\Lambda(j)b(m)a(mj)}{mj^{1+i\delta}}
     \right)
      + O \left( T L^{A+B+4} \log L \right) \ .  \\
   \label{eq:bunny}
\end{split}
\end{equation}
\end{s1}
\pr  By swapping summation order
\begin{equation}
  I  = \sum_{m \le \frac{T}{2 \pi}} \sum_{n \le \frac{T}{2 \pi}}
        \frac{a(m)b(n)}{n} \left( \frac{n}{m} \right)^{i \delta}
        \sum_{2 \pi \max(m,n) \le \gamma \le T} \left( \frac{n}{m}
        \right)^{\rho} \ .
\end{equation}
We decompose $I=I_{1}+I_{2}+I_{3}$ where
\begin{equation}
   I_{1} = \sum_{m \le \frac{T}{2 \pi}} \frac{a(m)b(m)}{m}(N(T)-N(2 \pi m)) \ ,
           \label{eq:horla}
\end{equation}
\begin{equation}
   I_{2} = \sum_{m \le \frac{T}{2 \pi}} \sum_{m < n}
        \frac{a(m)b(n)}{n} \left( \frac{n}{m} \right)^{i \delta}
        \sum_{2 \pi n \le \gamma \le T} \left( \frac{n}{m}
        \right)^{\rho} \ ,
   \label{eq:i2}
\end{equation}
and $I_{3}$ is the remaining piece consisting of terms with $n
<m$. The second expression in~(\ref{eq:horla}) is $ \ll L^{A+B+1}
\sum_{m \le T} d^{2}(m) \ll TL^{A+B+4}$ and since $N(T) =
\frac{TL}{2 \pi} + O(T)$ we deduce
\begin{equation}
   I_{1} = \frac{T}{2 \pi} \log \left( \frac{T}{2 \pi} \right)
   \sum_{m \le \frac{T}{2 \pi}} \frac{a(m)b(m)}{m}
   + O(TL^{A+B+4}) \ .
   \label{eq:auteur}
\end{equation}
Note that for $u \in \mathbb{R}$ and $0 < C < T$
\begin{equation}
   \overline{\sum_{C < \gamma < T} u^{\rho}}
   = u \sum_{C < \gamma < T} \left( \frac{1}{u} \right)^{\rho}
   \label{eq:vache}
\end{equation}
which follows from the symmetry of the zeros about $\mathrm{Re}(s) =
 \frac{1}{2}$.  Consequently, we deduce
\begin{equation}
  \overline{I_{3}}
   =  \sum_{m \le \frac{T}{2 \pi}} \sum_{n < m}
        \frac{a(m)b(n)}{m} \left( \frac{m}{n} \right)^{i \delta}
        \sum_{2 \pi m \le \gamma \le T} \left( \frac{m}{n}
        \right)^{\rho}   \ .
   \label{eq:i3}
\end{equation}
This expression has the same form as $I_{2}$ except the
roles of $a(n)$ and $b(n)$ have been switched.  Thus the evaluation of
$I_{3}$ follows along similar lines to $I_{2}$.
Putting $x=\frac{n}{m}$ and noticing $n \ll T$,~(\ref{eq:gl}) implies
\begin{equation}
\begin{split}
  \sum_{2 \pi n \le \gamma \le T} x^{\rho}
  & = \left(-\frac{T}{2 \pi} + n \right) \Lambda(x)
  +  O \left( x (\log(2xT)) ( \log \log 3x)
             \right) \\
 &  +O \left( (\log x)
    \min \left( T, \frac{x}{\langle x \rangle} \right)
    \right)
  +
    O \left( (\log (2T)) \min \left( T, \frac{1}{\log x} \right)
    \right) \ . \\
 \label{eq:vert}
\end{split}
\end{equation}
By inserting~(\ref{eq:vert}) into the inner sum of~(\ref{eq:i2})
we obtain $I_{2} = I_{21} + I_{22} + I_{23} + I_{24}$ where
\begin{equation}
\begin{split}
  I_{21} & =
        \sum_{m \le \frac{T}{2 \pi}}
       \sum_{{\begin{substack}{ m < n \le \frac{T}{2 \pi}
         \\ m \mid n }\end{substack}}}
      \frac{a(m)b(n)}{n} \left( \frac{n}{m} \right)^{i \delta}
      \left( -\frac{T}{2 \pi} + n  \right)
      \Lambda \left( \frac{n}{m} \right)  \\
\end{split}
\end{equation}
and $I_{22}$-$I_{24}$ correspond to the other terms
in~(\ref{eq:vert}). Applying $d(uv) \le d(u)d(v)$ the second part
of this expression is
\[
   \ll L^{A+B} \sum_{jm \le \frac{T}{2 \pi}}
     d(m)d(mj) \Lambda(j)
  \ll T L^{A+B+3}
    \sum_{j \le T} \frac{\Lambda(j) d(j)}{j}  \ll
  T L^{A+B+4}
\]
since the final sum is $\ll \sum_{p \le T} \frac{\log p}{p}$.
Consequently, we deduce that
\begin{equation}
\begin{split}
 I_{21} & = - \frac{T}{2 \pi}
      \sum_{mj \le \frac{T}{2 \pi}} \frac{\Lambda(j)a(m)b(mj)}{mj^{1-i\delta}}
      + O \left( T L^{A+B+4} \right) \ . \\
\end{split}
\end{equation}
The next term is
\begin{equation}
\begin{split}
  I_{22}
  & \ll   \sum_{n \le \frac{T}{2 \pi}} \sum_{m < n}
      \frac{a(m)b(n)}{n} \left(
      \frac{n}{m}  \log \log \left( \frac{n}{m} \right)
      \left( \log\left(2\frac{n}{m}T \right)
       \right)
      \right) \\
  &  \ll T L^{A+B+1} \log L
     \sum_{n \le \frac{T}{2 \pi}} d(n) \sum_{m < n} \frac{d(m)}{m}
     \ll T L^{A+B+4} \log L  . \\
\end{split}
\end{equation}
The third term, $I_{23}$, is bounded by
\begin{equation}
      \sum_{n \le \frac{T}{2 \pi}} \sum_{m < n}
      \frac{a(m)b(n)}{n}  \log \left( \frac{n}{m} \right)
      \left(
      \min \left( T, \frac{\frac{n}{m}}{ \langle
      \frac{n}{m} \rangle} \right)
       \right) \ll L^{A+B}
      \sum_{m < n \le \frac{T}{2 \pi}} \frac{d(m)}{m}
      \frac{d(n) \log \frac{n}{m}}{ \langle
      \frac{n}{m} \rangle} \ .
      \label{eq:mignon}
\end{equation}
In the last sum in~(\ref{eq:mignon}), pairs $(m,n)$ such that
$\langle \frac{n}{m} \rangle > \frac{1}{4}$ contribute
\[ \ll  L^{A+B}
      \left( \sum_{m < \frac{T}{2 \pi}} \frac{d(m)}{m} \right)
      \left( \sum_{n \le \frac{T}{2 \pi}} d(n) \log n \right)
   \ll T L^{A+B+4}  \ .
\]
The remaining pairs satisfy
$\langle \frac{n}{m} \rangle \le \frac{1}{4}$.
For each pair $(m,n)$ with $m < n$ we uniquely write
$n = qm + r$ with $-\frac{m}{2} < r \le  \frac{m}{2}$
and $q = \lfloor \frac{n}{m} \rfloor$ or
$q  = \lfloor \frac{n}{m} \rfloor +1 $.
By the identity
\[
  \left\langle \frac{n}{m} \right\rangle =
   \left\langle q + \frac{r}{m} \right\rangle =
   \left\{ \begin{array}{ccl}
           \frac{|r|}{m} & \mathrm{if} & q = p^{k} \ \mathrm{and} \
           r \ne 0  \\
           \ge \frac{1}{2} & \mathrm{if} & q \ne p^{k}
           \ \mathrm{or} \ r=0  \\
           \end{array}
   \right.
   \label{eq:coq}
\]
we need only consider $q$ a prime power.
Thus the contribution from those terms with
$\langle \frac{n}{m} \rangle \le \frac{1}{4}$
in the final sum in~(\ref{eq:mignon}) is
\begin{equation}
   \sigma := \sum_{m < X} d(m)
   \sum_{q \le \frac{X}{m} +1, \ q=p^{k}} \log(q)
   \sum_{1 \le |r| \le \frac{m}{4}} \frac{d(qm+r)}{|r|}
\end{equation}
where $X =  \frac{T}{2 \pi}$
and hence $I_{23} \ll \sigma L^{A+B} + TL^{A+B+4}$.  Furthermore,
we write $\sigma = \sigma_{1}+\sigma_{2}+\sigma_{3}$
according to the cases:
$(i) \ q = p \ \mathrm{prime}$,
$(ii) \ q = p^{k}, \ k \ge 2, \ X^{1-\delta} < p^{k} \le X+1$, and
$(iii) \ q = p^{k}, \ k \ge 2, \ p^{k} \le X^{1-\delta}$ where
$\delta$ is a small positive constant.
The contribution from $(i)$ is
\begin{equation}
\begin{split}
   \sigma_{1} & \ll \sum_{m < X} d(m)
                  \sum_{q \le \frac{X}{m}+1} \Lambda(q)
                  \sum_{1 \le |r| \le \frac{m}{4}} \frac{d(qm+r)}{|r|}
              \\
   &  \ll \sum_{n \le X} d(n)
                      \sum_{|r| \le \frac{X}{4}, |r| < n} \frac{1}{|r|}
                      \sum_{qm = n-r} \Lambda(q) d(m)  \\
\end{split}
\end{equation}
where we wrote $n=qm + r$ and noticed that $qm \le 2X$.
Since $d(m) \le d(qm) = d(n-r)$ and
$\sum_{qm = n-r} \Lambda(q) = \log(n-r) \le \log(2X) \ll L$,
our sum is bounded by
\begin{equation}
                \ll L \sum_{|r| \le \frac{X}{4}} \frac{1}{|r|}
                      \sum_{|r| < n \le X}
                      d(n)d(n-r)
                  \ll T L^{3}
                  \sum_{r \le \frac{X}{4}}
                  \frac{\sigma_{-1}(r)}{r}  \ .
   \label{eq:e38}
\end{equation}
The right-most inequality follows by Lemma 2 and thus
\begin{equation}
    \sigma_{1} \ll T L^3 \sum_{r \le \frac{X}{4}} \frac{1}{r}
                    \sum_{g \mid r} \frac{1}{g}
                = T L^3 \sum_{g \le \frac{X}{4}} \frac{1}{g^2}
                    \sum_{s \le \frac{X}{4}} \frac{1}{s}
                  \ll T L^4 \ .
   \label{eq:e39}
\end{equation}
Observe that in $\sigma_{2}$, condition $(ii)$ implies $ m \ll
\frac{X}{q} \le X^{\delta}$ and since $d(qm+r) \ll X^{\delta}$ we
have
\begin{equation}
\begin{split}
  \sigma_{2}
  & \ll \sum_{m \le X^{\delta}} d(m) \sum_{p^{k} \ll \frac{X}{m}, \ k
  \ge 2} \log(p^{k}) \sum_{r \le \frac{m}{4}} \frac{d(qm+r)}{r} \\
  & \ll X^{\delta} \log X
   \sum_{m \le X^{\delta}} d(m) \sum_{p^{k} \ll \frac{X}{m}, \ k
  \ge 2} \log(p^{k}) \ll T^{\frac{1}{2}+\delta} \ . \\
\end{split}
\end{equation}
In the final piece we have
\begin{equation}
\begin{split}
  & \sigma_{3} \ll  \sum_{p^{k} \le X^{1-\delta}, \ k \ge 2} \log(p^{k})
  \sum_{|r| \le \frac{X}{2}} \frac{1}{|r|}
   \sum_{{\begin{substack}{2|r| \le m < X
         \\ m \ll \frac{X}{p^{k}} \ , \ p^{k}m+r < X}\end{substack}}}
  d(m) d(p^{k}m+r) \ .  \\
  \label{eq:ours}
\end{split}
\end{equation}
By Cauchy-Schwarz, the inner sum in~(\ref{eq:ours}) is
\begin{equation}
  \ll
  \left( \frac{X}{p^{k}} \right)^{\frac{1}{2}}
   \log^{\frac{3}{2}} X
  \left(
  \sum_{{\begin{substack}{n \ll X
         \\ n \equiv r \ \mathrm{mod} \ p^{k}}\end{substack}}}
   d(n)^{2}
  \right)^{\frac{1}{2}} \ . \\
  \label{eq:chat}
\end{equation}
We now establish
\begin{equation}
    \sum_{{\begin{substack}{n \ll X
         \\ n \equiv r \ \mathrm{mod} \ p^{k}}\end{substack}}}
  d(n)^{2}
  \ll \min( \mathrm{ord}_{p}(r),k)^{2} \
  \frac{X (\log X)^{3}}{p^{k}} \ .
  \label{eq:faim}
\end{equation}
If $(r,p) = 1$~(\ref{eq:faim}) is true by~(\ref{eq:vino}).
On the other hand, suppose $(r,p) > 1$ and $r = p^{u}s$ with $(s,p) =1$.
If $u \ge k$ then we have
\[
    \sum_{{\begin{substack}{n \ll X
         \\ n \equiv 0 \ \mathrm{mod} \ p^{k}}\end{substack}}}
   d(n)^{2}
  \le d(p^{k})^{2}  \sum_{j \ll \frac{X}{p^{k}}} d(j)^{2}
  \ll \frac{k^{2}}{p^{k}} (X  \log^{3} X) \ .
\]
In the case $1 \le u < k$, an analogous calculation establishes
the other bound in~(\ref{eq:faim}).
Combining~(\ref{eq:ours}),~(\ref{eq:chat}), and~(\ref{eq:faim})
we have
\[
  \sigma_{3} \ll
  (X \log^{4} X) \sum_{p^{k} \le X^{1-\delta}, \ k \ge 2}
  \frac{k \log(p^{k})}{p^{k}} \ll TL^{4} \ .
\]
Putting together our estimates for the $\sigma_{i}$, we have
$\sigma \ll TL^{4}$ and hence $I_{23} \ll TL^{A+B+4}$. Finally,
$I_{24}$ is
\[
  \ll L^{1+A+B} \sum_{n \le \frac{T}{2 \pi}} \sum_{m < n}
      \frac{d(m)d(n)}{n} \frac{1}{\log \frac{n}{m}}
      \ll  L^{1+A+B} \sum_{r < \frac{T}{2 \pi}} \frac{1}{r}
        \sum_{r < n \le \frac{T}{2 \pi}}
        d(n-r)d(n) \ .
\]
Notice that the last sum was already treated in~(\ref{eq:e38})
and~(\ref{eq:e39}), so we have $I_{24} \ll T L^{A+B+4}$.
Thus we arrive at
\begin{equation}
   I_{2} = - \frac{T}{2 \pi}
      \sum_{mj \le \frac{T}{2 \pi}} \frac{\Lambda(j)a(m)b(mj)}{mj^{1-i\delta}}
      + O \left( T L^{A+B+4} \log L \right) \ .
      \label{eq:bonne}
\end{equation}
Starting from~(\ref{eq:i3}) an analogous calculation demonstrates that
\begin{equation}
   I_{3} = - \frac{T}{2 \pi}
      \sum_{mj \le \frac{T}{2 \pi}} \frac{\Lambda(j)b(m)a(mj)}{mj^{1+i\delta}}
      + O \left( T L^{A+B+4} \log L \right) \ .
      \label{eq:corps}
\end{equation}
Combining~(\ref{eq:auteur}),~(\ref{eq:bonne}), and~(\ref{eq:corps}) finishes the proof of the lemma. \\

In the next lemma, we evaluate the
second and third sums of~(\ref{eq:bunny}).
\newtheorem{const}[gonek]{Lemma}
\begin{const}
Suppose we have two sequences
$a(n) \ll \log^{A} (n) d(n)$ and $b(n) \ll \log^{B} (n)  d(n)$
which satisfy for each $p \le t$
\begin{equation}
   \sum_{n \le t} \frac{a(n)b(pn)}{n} =
   \sum_{u+v= \beta} s_{uv} \log^{u} p \log^{v} t
   + O \left( \log^{\beta - 1} t
              + \frac{(\log^C p) (\log^{\beta} t)}{p} \right)
  \label{eq:prod}
\end{equation}
where  $\beta,C$ are positive absolute constants, $u,v \ge 0$,
$s_{uv} \in \mathbb{C}$, and the implied constant in the error
term depends only on $a(n)$ and $b(n)$. We associate to an
expansion of the form~(\ref{eq:prod}) the constant
\begin{equation}
   \mathcal{A}(a,b) = \sum_{u+v=\beta} s_{uv}
   \frac{u! v!}{(u+v+1)!}
   \ .
   \label{eq:Aab}
\end{equation}
Then we have
\begin{equation}
\begin{split}
   & M(a,b;X,\delta)
   := \sum_{mk \le X} \frac{\Lambda(k) a(m) b(mk)}{k^{1-i\delta} m} \\
  &  = \tilde{L}^{\beta+1}
  \sum_{k=0}^{\infty} \frac{(i \delta \tilde{L})^{k}}{k!}
     \sum_{u+v=\beta} s_{uv} \frac{(u+k)! v!}{(u+v+k+1)!}
  + O(\tilde{L}^{\max(\beta,A+B+4)})  \\
\end{split}
\end{equation}
where $\tilde{L} = \log X$ and $\delta \in \mathbb{R}$.
Moreover if $\delta=0$, this reduces to
\begin{equation}
    M(a,b;X,0) = \sum_{mk \le X} \frac{\Lambda(k) a(m) b(mk)}{k m}
   = \tilde{L}^{\beta+1}   \mathcal{A}(a,b)
  + O(\tilde{L}^{\max(\beta,A+B+4)}) \ .
\end{equation}
\end{const}
\pr
In the sum $M(a,b;X,\delta)$ the prime powers
$p^{\alpha}$ with $\alpha \ge 2$ contribute
\[
  \sum_{p^{\alpha} \le X, \alpha \ge 2}
  \frac{\Lambda(p^{\alpha})}{p^{\alpha}}
   \sum_{m < \frac{X}{p^{\alpha}}} \frac{a(m)b(mp^{\alpha})}{m} \ll
  \tilde{L}^{A+B+4}\sum_{p^{\alpha} \le X ,\ \alpha \ge 2}
  \frac{\alpha  \Lambda(p^{\alpha})}{p^{\alpha}}
  \ll \tilde{L}^{A+B+4} \ .
\]
We arrive at
\[ M(a,b;X,\delta) =
      \sum_{p \le X} \frac{\Lambda(p)}{p^{1-i\delta}}
      \sum_{m < \frac{X}{p}} \frac{a(m)b(mp)}{m}
      + O(\tilde{L}^{A+B+4}) \ .
\]
We replace the inner sum above by the expression on the right side
of~(\ref{eq:prod}). The contribution to $M(a,b;X,\delta)$ coming
from the error term in~(\ref{eq:prod}) is
\[ \tilde{L}^{\beta-1} \sum_{p \le X} \frac{\log p}{p}
   + \tilde{L}^{\beta} \sum_{p \ge 2} \frac{\Lambda(p) \log^C p}{p^2}
   \ll \tilde{L}^{\beta}  \ .
\]
This demonstrates that
\begin{equation}
  M(a,b;X,\delta) =
  \sum_{u+v=\beta} s_{uv} \sum_{p \le X} \frac{\Lambda(p)}{p^{1-i\delta}}
    \log^{u} p \ \log^{v} \left( \frac{X}{p} \right)
  + O(\tilde{L}^{\max{\beta,A+B+4}}) \ .
  \label{eq:vers}
\end{equation}
By Stieltjes integration,
\begin{equation}
  \sum_{p \le X} \frac{\Lambda(p)}{p^{1-i\delta}}
    \log^{u} p \ \log^{v} \left( \frac{X}{p} \right)
  = \int_{1}^{X} \log^{u} t
      \ \log^{v} \left( \frac{X}{t} \right)
      \ \frac{d \theta(t)}{t^{1-i\delta}}
  \label{eq:dij}
\end{equation}
where $\theta(t) = \sum_{p \le t} \log p$ .
The prime number theorem is $\theta(t) = t + O(t\exp(-c \sqrt{\log
  t}))$ and thus the main part of~(\ref{eq:dij}) equals
\begin{equation}
    \int_{1}^{X} \log^{u} t
      \ \log^{v} \left( \frac{X}{t} \right)
      \ \frac{dt}{t^{1-i \delta}}
  =
    \sum_{k=0}^{\infty} \frac{(i \delta)^{k}}{k!}
    \int_{1}^{X} \log^{u+k} t
      \ \log^{v} \left( \frac{X}{t} \right)
      \ \frac{dt}{t}
    \label{eq:fourmi}
\end{equation}
\[
   = \tilde{L}^{\beta+1}  \sum_{k=0}^{\infty}
      \frac{(i \delta \tilde{L})^{k}}{k!}
      \int_{0}^{1} x^{u+k}(1 - x)^{v} \ dx
    = \tilde{L}^{\beta+1}  \sum_{k=0}^{\infty}
      \frac{(i \delta \tilde{L})^{k}}{k!}
      \frac{(u+k)! v!}{(u+v+k+1)!}
\]
where we made the variable change $x = (\log t)/ \tilde{L}$. The
contribution arising from the error term in the prime number
theorem is  easily seen to be $\tilde{L}^{\beta}$.
Combining~(\ref{eq:vers}),~(\ref{eq:dij}), and~(\ref{eq:fourmi})
establishes the lemma.  \\

Putting together Lemmas 3 and 4 we have the following computation
of the main term of $I(a,b;T)$ in~(\ref{eq:iabt}) subject to
various conditions on the sequences $a(n)$ and $b(n)$.
\newtheorem{comb}[gonek]{Lemma}
\begin{comb}
Suppose we have
two sequences $a(n) \ll \log^{A}(n) d(n)$
and $b(n) \ll \log^{B}(n) d(n)$ such that
\begin{equation}
  \sum_{n \le t} \frac{a(n)b(n)}{n} = c_{a,b} \log^{\beta} t +
                                        O(\log^{\beta-1} t) \ ,
\end{equation}
\begin{equation}
   \sum_{n \le t} \frac{a(n)b(pn)}{n} =
   \sum_{u+v= \beta} s_{uv} \log^{u} p \log^{v} t
   + O \left( \log^{\beta - 1} t
              + \frac{(\log^C p) ( \log^{\beta} t)}{p} \right) \ , \
              and \
  \label{eq:abpn}
\end{equation}
\begin{equation}
   \sum_{n \le t} \frac{b(n)a(pn)}{n} =
   \sum_{u+v= \beta} t_{uv} \log^{u} p \log^{v} t
   + O \left( \log^{\beta - 1} t
              + \frac{(\log^C p) (\log^{\beta} t)}{p} \right) \ ,
  \label{eq:abpn2}
\end{equation}
where $c_{a,b},A,B,\beta,C$ are fixed positive constants.
Moreover, suppose that~(\ref{eq:abpn}) and~(\ref{eq:abpn2}) hold
for $p \le t$ and the constant in the error term is independent of
$p$. Then we have
\begin{equation}
   I(a,b;T)  = \frac{TL^{\beta+1}}{2 \pi}
              \left( c_{a,b} - \mathcal{A}(a,b) - \mathcal{A}(b,a) \right)
              + O(T(L^{\max(\beta,A+B+4)}))
\end{equation}
where $\mathcal{A}(a,b)$ and $\mathcal{A}(b,a)$ are constants
defined by~(\ref{eq:Aab}).
\end{comb}
{\bf More notation} For arbitrary sequences
$a(n)$ and $b(n)$ define the functions
\begin{equation}
   T_{a,b}(t) = \sum_{n \le t} \frac{a(n)b(n)}{n} \ \mathrm{and} \
   T_{a,b;p}(t) = \sum_{n \le t} \frac{a(n)b(pn)}{n} \ .
   \label{eq:tabt}
\end{equation}
Furthermore, we use the simplified notation
\begin{equation}
    T_{\mu,\nu}(t) := T_{d^{(\mu)},d^{(\nu)}}(t) \ \mathrm{and} \
    T_{\mu,\nu;p}(t) := T_{d^{(\mu)},d^{(\nu)};p}(t)
    \label{eq:tmnp}
\end{equation}
for $\mu, \nu \in \mathbb{Z}_{\ge 0}$. Also define
\begin{equation}
    T_{(n_{1},n_{2}),(n_{3},n_{4})}(t)
    := T_{d^{(n_{1},n_{2})},d^{(n_{3},n_{4})}(n) }
    = \sum_{n \le t} \frac{
    d^{(n_{1},n_{2})}(n)d^{(n_{3},n_{4})}(n)}{n}
    \label{eq:t1234}
\end{equation}
for $n_{1},n_{2},n_{3},n_{4} \in \mathbb{Z}_{\ge 0}$. Recall that
$d^{(\mu,\nu)}(n)$ is defined by~(\ref{eq:duvc}) and $d^{(\mu)}(n)
= d^{(\mu,0)}(n)$. Note that $T_{(\mu,0),(\nu,0)}(t) =
T_{\mu,\nu}(t)$. By Lemma 5, we need to evaluate sums of the
form~(\ref{eq:tabt}) in order to compute the constants
$\mathcal{A}(a,b)$ in~(\ref{eq:Aab}). Once this is done we obtain
the
main term asymptotic for $I(a,b;T)$ in~(\ref{eq:iabt}). \\

Our calculations require an effective version of
Perron's formula.
\newtheorem{per}[gonek]{Lemma}
\begin{per}
Let $F(s) := \sum_{n\ge 1} a_{n} n^{-s}$ be a Dirichlet series with
finite abscissa of absolute convergence $\sigma_{a}$.  Suppose
there exists a real number $\alpha \ge 0$ such that
\begin{equation}
  \sum_{n=1}^{\infty} |a_{n}|n^{-\sigma} \ll (\sigma-\sigma_{a})^{-\alpha}
  \
  (\sigma > \sigma_{a})
\end{equation}
and that B is a non-decreasing function such that $|a_{n}| \le B(n)$
for $n \ge 1$.  Then for $x \ge 2, T \ge 2, \sigma \le \sigma_{a},
\kappa := \sigma_{a}-\sigma+(\log x)^{-1}$, we have
\begin{equation}
\begin{split}
  \sum_{n \le x} \frac{a_{n}}{n^{s}}
  & = \frac{1}{2 \pi i} \int_{\kappa-iT}^{\kappa+iT} F(s+w)
  \frac{x^{w}}{w} \, dw
  + O \left( \frac{x^{\sigma_{a}-\sigma}(\log x)^{\alpha}}{T}
             + \frac{B(2x)}{x^{\sigma}}
             \left( 1 + x \frac{\log T}{T} \right)
  \right) \ .  \\
\end{split}
\end{equation}
\end{per}
\pr This is Corollary 2.1 p.133 of \cite{Te}. \\

The evaluation of~(\ref{eq:t1234}) follows closely Theorem 7 of
\cite{Ha} pp.296-297.
\newtheorem{me}[gonek]{Lemma}
\begin{me}
We have $ T_{(n_{1},n_{2}),(n_{3},n_{4})}(t)
   = P(\log t) + O_{\epsilon}(t^{-\frac{1}{2}+\epsilon})$
where $P(x)$ is a polynomial of degree $n_{1}+n_{2}+n_{3}+n_{4}+4$
with leading coefficient
\begin{equation}
   \frac{a_{2}n_{1}!n_{2}!n_{3}!n_{4}!}{(n_{1}+n_{2}+n_{3}+n_{4}+4)!}
   \sum_{a=0}^{n_{1}} \sum_{c=0}^{n_{2}}
   \binom{n_{3}+1+a+c}{n_{3}}
   \binom{n_{4}+1+n_{1}+n_{2}-a-c}{n_{4}} \ .
   \label{eq:comb}
\end{equation}
A special case of this result is
$T_{\mu,\nu}(t) =
   Q(\log t) + O_{\epsilon}(t^{-\frac{1}{2}+\epsilon})$
where $Q(x)$ is a polynomial of degree $\mu+\nu+4$ with leading
coefficient
\begin{equation}
   C(\mu,\nu) = \frac{\mu! \nu!}{(\mu+\nu+4)!}
  \left( \binom{\mu+\nu+2}{\mu+1} - 1 \right) \ .
  \label{eq:cuk}
\end{equation}
\end{me}
\pr Define $\sigma_{u,v}(n) = \sum_{d_{1}d_{2}=n} d_{1}^{u}
d_{2}^{v}$ where $u,v \in \mathbb{C}$. Let
$z_{1},z_{2},z_{3},z_{4} \in \mathbb{C}$ and define the Dirichlet
series
\[
   F(s;\vec{z}) := \sum_{n \le t} \frac{ \sigma_{-z_{1},-z_{2}}(n)
                 \sigma_{-z_{3},-z_{4}}(n)}{n^{s+1}}
\]
where $\vec{z} = (z_{1},z_{2},z_{3},z_{4})$.
Observe the relationship
\begin{equation}
   \left.
   (-1)^{n_{1}+n_{2}+n_{3}+n_{4}}
   \frac{d^{n_{1}}}{dz_{1}^{n_{1}}}
   \frac{d^{n_{2}}}{dz_{2}^{n_{2}}}
   \frac{d^{n_{3}}}{dz_{3}^{n_{3}}}
   \frac{d^{n_{4}}}{dz_{4}^{n_{4}}} F(s;\vec{z}) \right|_{\vec{z} =
   \vec{0}} =
   \sum_{n \ge 1}
        \frac{ d^{(n_{1},n_{2})}(n)d^{(n_{3},n_{4})}(n)}{n^{s+1}}
   \ .
   \label{eq:idn}
\end{equation}
We denote the generating function in~(\ref{eq:idn}) $F(s)$.
On the other hand, by Ramanujan's calculation (see \cite{Ti} pp.8-9),
$F(s;\vec{z})$ equals
\begin{equation}
      \frac{ \zeta(1+s+z_{2}+z_{4}) \zeta(1+s+z_{1}+z_{4})
             \zeta(1+s+z_{2}+z_{3})  \zeta(1+s+z_{1}+z_{3}) }
      {\zeta(2+2s+z_{1}+z_{2} +z_{3}+z_{4})} \ .
      \label{eq:ram}
\end{equation}
By~(\ref{eq:idn}) and~(\ref{eq:ram}) we deduce that
\begin{equation}
  F(s) := (-1)^{N}
  \sum_{\vec{a} \in (\mathbb{Z}_{\ge 0})^{4}}
  G_{\vec{a}}(s)
  \zeta^{(a_{1})}(1+s) \zeta^{(a_{2})}(1+s)
  \zeta^{(a_{3})}(1+s) \zeta^{(a_{4})}(1+s)
  \label{eq:idn2}
\end{equation}
where $N = n_{1}+n_{2}+n_{3}+n_{4}$ and
$\vec{a} =(a_{1},a_{2},a_{3},a_{4}) \in (\mathbb{Z}_{\ge 0})^{4}$
ranges over a finite sum.
Moreover, the functions
$G_{\vec{a}}(s)$ have absolutely convergent Dirichlet series
in $\mathrm{Re}(s) > -\frac{1}{2}$.
A careful examination of~(\ref{eq:ram})
reveals that the leading term in the Laurent expansion of
$F(s)$ derives from the expression
\begin{equation}
  (-1)^{N}
   \frac{d^{n_{1}}}{dz_{1}^{n_{1}}}
   \frac{d^{n_{2}}}{dz_{2}^{n_{2}}}
   \frac{d^{n_{3}}}{dz_{3}^{n_{3}}}
   \frac{d^{n_{4}}}{dz_{4}^{n_{4}}} G(1+s;\vec{z})
   \label{eq:ferme}
\end{equation}
where
\[
   G(w;\vec{z}) :=
   \zeta(w+z_{2}+z_{4}) \zeta(w+z_{1}+z_{4})
             \zeta(w+z_{2}+z_{3})  \zeta(w+z_{1}+z_{3})  \ .
\]
An application of the product rule
$(f(z)g(z))^{(n)} := \sum_{j=0}^{n} \binom{n}{j}
  f^{(j)}(z) g^{(n-j)}(z)$
in each of the variables to~(\ref{eq:ferme}) yields
\begin{equation}
\begin{split}
  & (-1)^{N} \sum_{a,c,e,g} \binom{n_{1}}{a}
  \binom{n_{2}}{c}
  \binom{n_{3}}{e}
  \binom{n_{4}}{g}
  \zeta^{(a+e)}(1+s+z_{1}+z_{3})
   \cdot \\
  & \cdot \zeta^{(c+f)}(1+s+z_{2}+z_{3})
    \zeta^{(b+g)}(1+s+z_{1}+z_{4})
    \zeta^{(d+h)}(1+s+z_{2}+z_{4}) \\
\end{split}
\end{equation}
where $a+b=n_{1},c+d=n_{2},e+f=n_{3}$, and $g+h=n_{4}$. Thus
\begin{equation}
\begin{split}
  F(s) := & \frac{(-1)^{N}}{\zeta(2+2s)}
   \sum_{a,c,e,g} \binom{n_{1}}{a}
   \binom{n_{2}}{c}
   \binom{n_{3}}{e}
   \binom{n_{4}}{g} \cdot \\
  & \zeta^{(a+e)}(1+s)
   \zeta^{(c+f)}(1+s)
    \zeta^{(b+g)}(1+s)
    \zeta^{(d+h)}(1+s) + R(s)
  \label{eq:la0} \\
\end{split}
\end{equation}
where $R(s)$ is a function with a pole of order at most $N+3$ at
$s=0$. Note that we have the expansions $  \zeta^{(k)}(1+s) =
\frac{(-1)^{k}k!}{s^{k+1}} + c_{k} + \cdots$ and
$\frac{1}{\zeta(2+2s)} = \frac{1}{\zeta(2)} + c^{'}s + \cdots$ for
constants $c_{k}$ and $c^{'}$. By combining~(\ref{eq:la0}) with
these expansions, we have $F(s) = \frac{6}{\pi^{2}}Cs^{-N-4} +
R^{'}(s)$. Here $R^{'}(s)$ consists of those terms in the Laurent
expansion with pole of order at most $N+3$ and
\begin{equation}
   C =
      \sum_{a,c,e,g}
  \binom{n_{1}}{a} \binom{n_{2}}{c} \binom{n_{3}}{e}
          \binom{n_{4}}{g} \cdot
  (a+e)! (c+f)! (b+g)! (d+h)!
  \label{eq:lapin}
\end{equation}
We simplify $C$ by applying the identity
\begin{equation}
   \sum_{k=0}^{l} \binom{l-k}{m} \binom{q+k}{n}
   = \binom{l+q+1}{m+n+1}
   \label{eq:biid}
\end{equation}
valid for integers $l,m \ge 0$ and integers $n \ge q \ge 0$
(see \cite{GKP} p.169). The sum over $g$ in~(\ref{eq:lapin}) is
\begin{equation}
\begin{split}
   & \sum_{g=0}^{n_{4}} \binom{n_{4}}{g} (b+g)! (d+h)! =
   \sum_{g=0}^{n_{4}} \frac{n_{4}!}{g!(n_{4}-g)!}
   (b+g)!(d+n_{4}-g)! \\
   & = n_{4}!b!d! \sum_{g=0}^{d+n_{4}} \binom{b+g}{b}
   \binom{d+n_{4}-g}{d}
   = n_{4}!b!d! \binom{n_{4}+b+d+1}{b+d+1}  \\
\end{split}
\end{equation}
where we applied~(\ref{eq:biid}). Similarly, the sum over $e$ is
\begin{equation}
  \sum_{e=0}^{n_{3}} \binom{n_{3}}{e}(a+e)!(c+f)! =
  n_{3}!a!c! \binom{n_{3}+a+c+1}{a+c+1}  \ .
\end{equation}
Since $\binom{n_{1}}{a}a!b!= n_{1}!$ and
$\binom{n_{2}}{c}c!d!= n_{2}!$ the total sum is
\begin{equation}
  C  = n_{1}!n_{2}! n_{3}! n_{4}!
   \sum_{a=0}^{n_{1}} \sum_{c=0}^{n_{2}}
   \binom{n_{3}+a+c+1}{n_{3}}  \binom{n_{4}+1+n_{1}+n_{2}-a-c}
    {n_{4}} \ .
   \label{eq:con}
\end{equation}
This shows that $F(s)s^{-1} := \frac{C}{s^{N+5}} +
\frac{C^{'}}{s^{N+4}} + \cdots$ for constants $C$ and $C^{'}$.
Hence the residue of $F(s)t^{s}s^{-1}$ at $s=0$ is $P(\log t)$
where $P(t)$ is a polynomial of degree $N+4$ with leading
coefficient $\frac{6}{\pi^{2}} C/(N+4)!$. By Lemma 6 applied with
$\alpha = N+4$, $s=\sigma_{a}=1$, and $B(t) \ll_{\epsilon}
t^{\epsilon}$ it follows that
\begin{equation}
   T(t) := \frac{1}{2 \pi i}
        \int_{\kappa-iT}^{\kappa+iT} F(w) \frac{t^{w}}{w} \, dw
        + O \left(  \frac{(\log t)^{N+4}}{T} +
        \frac{1}{t^{1-\epsilon}}
        \left( 1 + t \frac{\log T}{T}         \right) \right)
\end{equation}
where $\kappa = (\log t)^{-1}$. By the residue theorem, the integral is
\begin{equation}
 P(\log t)
  -  \frac{1}{2 \pi i}
  \left( \int_{\kappa+iT}^{c+iT} + \int_{c+iT}^{c-iT} +
         \int_{c-iT}^{\kappa-iT} \right) \frac{F(w)t^{w}}{w} \, dw \ .
\end{equation}
where $c= -\frac{1}{2} + \epsilon$.
We only sketch how to estimate these integrals since the argument
is standard.
The first and third integral may be computed by
using known bounds for $\zeta(s)$ in the critical strip.
The second integral requires the result
\begin{equation}
   \int_{1}^{U} |\zeta^{(a)}(\tau + it)|^{4} \ll_{a,\tau} U
   \label{eq:mo4a}
\end{equation}
for $a \in \mathbb{Z}_{\ge 0}$ and $\tau > \frac{1}{2}$. This may
be proven by following the argument of Theorem 7.5 pp.146-147 of
\cite{Ti}. An appropriate choice of $T$ then yields an error term
of $t^{-\frac{1}{2}+\epsilon}$ to complete the proof.  For the
special case $T_{\mu,\nu}(t)$, we set $n_{1} = \mu, n_{2}=0, n_{3}
= \nu$, and $n_{4}=0$.  Applying the binomial identity (see
\cite{GKP} p.174)
\begin{equation}
   \sum_{k=0}^{n} \binom{r+k}{k} = \binom{r+n+1}{n}
   \label{eq:lampe}
\end{equation}
for $r,n \in \mathbb{Z}_{\ge 0}$,~(\ref{eq:comb}) reduces to
\begin{equation}
   \frac{6}{\pi^{2}} \frac{\mu! \nu!}{(\mu+\nu+4)!}
  \sum_{a=0}^{\mu} \binom{\nu+1+a}{\nu}
   = \frac{6}{\pi^{2}} \frac{\mu! \nu!}{(\mu+\nu+4)!}
  \left( \binom{\mu+\nu+2}{\mu+1} - 1 \right)
\end{equation}
and thus~(\ref{eq:cuk}) is verified. \\

We now record the special cases of Lemma 7 which are required in
the proof of Corollary 1. In Table 1, we associate to each pair of
sequences $(a,b)$ the main term of $T_{a,b}(t)$
in~(\ref{eq:tabt}).
\begin{center}
{\bf Table 1}
\end{center}
\begin{center}
\begin{tabular}{|c|c|c|c|}
\hline
 $(a,b)$ & $T_{a,b}(t)$
 & $(a,b)$ & $T_{a,b}(t)$ \\
\hline
\hline
$ (d,d) $ & $ \sim a_{2} \cdot
          \frac{1}{24} l^{4}$
  &  $ (d^{(1)},d^{(2)}) $ & $ \sim a_{2} \cdot
          \frac{1}{280}  l^{7}$ \\
\hline
$ (d,d^{(1)}) $ & $ \sim a_{2} \cdot
          \frac{1}{60}   l^{5}$
  & $ (d^{(1)},\alpha) $ & $ \sim a_{2}  \cdot
          \frac{1}{420}   l^{7} $ \\
\hline
$ (d^{(1)},d^{(1)}) $ & $ \sim a_{2} \cdot
          \frac{1}{144}   l^{6}$ & $ (d^{(2)},d^{(2)}) $ & $ \sim a_{2}  \cdot
          \frac{19}{10080}  l^{8}$ \\
\hline
$ (d,d^{(2)}) $ & $ \sim a_{2} \cdot
          \frac{1}{120}  l^{6}$  & $ (\alpha,\alpha) $ & $ \sim a_{2}  \cdot
          \frac{17}{20160}   l^{8}$ \\
\hline
$ (d,\alpha) $ & $ \sim a_{2} \cdot
          \frac{1}{180}   l^{6}$ & & \\
\hline
\end{tabular}
\end{center}

We now evaluate $T_{d^{(\mu)},d^{(\nu)};p}(t)$ in~(\ref{eq:tmnp}).
\newtheorem{Tuvp}[gonek]{Lemma}
\begin{Tuvp}
Let $\mu, \nu \ge 0$ be integers and let $p$ be any prime $\le t$.
We have
\begin{equation}
\begin{split}
  T_{\mu,\nu;p}(t) & = a_{2} \cdot
  \left(  2C(\mu,\nu) l^{\mu+\nu+4}
  + \sum_{k=0}^{\nu-1} \binom{\nu}{k} (\log^{\nu-k} p) l^{\mu+k+4}
    C(\mu,k) \right) \\
  & + O_{\mu,\nu}
  \left( \frac{l^{\mu+\nu+4}}{p} + l^{\mu+\nu+3} \right)
  \label{eq:tuvp}
\end{split}
\end{equation}
where $l = \log t$, $C(\mu,k)$ is defined by~(\ref{eq:cuk}), and
the sum only occurs if $\nu \ge 1$ and is zero otherwise.
\end{Tuvp}
\pr  First note that
\begin{equation}
    T_{\mu,\nu;p}(t) = \left . \frac{d^{\nu}}{dz^{\nu}}
   \sum_{n \le t} \frac{ d^{(\mu)}(n) \sigma_{z}(pn)}{n}
   \, \right |_{z=0}  \ .
   \label{eq:e1}
\end{equation}
Inserting the identity
\[ \sigma_{z}(n_{1} n_{2})  =
   \sum_{m \mid (n_{1},n_{2}) } \mu(m) m^{z}
   \sigma_{z} \left(\frac{n_{1}}{m} \right)
   \sigma_{z} \left( \frac{n_{2}}{m} \right)
\]
in~(\ref{eq:e1}) and inverting summations we obtain
\begin{equation}
   \sum_{n \le t} \frac{ d^{(\mu)}(n) \sigma_{z}(pn)}{n}
   = \sum_{j \le t}
      \frac{d^{(\mu)}(j)\sigma_{z}(p)\sigma_{z}(j)}{j} - \frac{1}{p}
      \sum_{j \le \frac{t}{p}}
      \frac{d^{(\mu)}(pj)p^{z}\sigma_{z}(j)}{j} \ .
   \label{eq:pain}
\end{equation}
Observe that
\begin{equation}
  \left . \frac{d^{\nu}}{dz^{\nu}}  \sigma_{z}(p)\sigma_{z}(j) \right|_{z=0}
   = \sum_{k=0}^{\nu-1} \binom{\nu}{k} \log^{\nu-k}(p)
  d^{(k)}(j) + 2 d^{(\nu)}(j)  \ and \
  \label{eq:e2}
\end{equation}
\begin{equation}
   \left . \frac{d^{\nu}}{dz^{\nu}}  p^{z} \sigma_{z}(j) \right|_{z=0}
   = \sum_{k=0}^{\nu} \binom{\nu}{k} \log^{\nu-k}(p)
  d^{(k)}(j)
  \label{eq:e3}
\end{equation}
where~(\ref{eq:e3}) follows from~(\ref{eq:e2}) since
$\sigma_{z}(p) = 1 + p^{z}$.
Combining~(\ref{eq:e1}),~(\ref{eq:pain}),~(\ref{eq:e2}),
and~(\ref{eq:e3}) we deduce that
 $T_{\mu,\nu;p}(t)$ equals
\begin{equation}
 2 T_{\mu, \nu}(t)
  + \sum_{k=0}^{\nu-1} \binom{\nu}{k} (\log^{\nu-k} p) \
  T_{\mu,k}(t)
  - \frac{1}{p} \left(
  \sum_{k=0}^{\nu} \binom{\nu}{k} (\log^{\nu-k} p) \, T_{k,\mu;p}
  \left( \frac{t}{p} \right) \right)  \ .
  \label{eq:stylo}
\end{equation}
The trivial bound $ T_{k,\mu;p}(t) \ll (\log^{k+4} t) ( \log^{\mu}
(pt))$ follows from $d^{(k)}(j) \le (\log^{k} j) d(j)$ and hence
the error term is $\ll p^{-1} \log^{\nu+\mu+4} t$. Applying Lemma
7 to each expression in the main term of~(\ref{eq:stylo}) completes the proof of the lemma. \\

We now compute $\mathcal{A}(d^{(\mu)},d^{(\nu)})$.
\newtheorem{auv}[gonek]{Lemma}
\begin{auv}
Let $\mu,\nu \ge 0$ be integers.  We have
\begin{equation}
   \mathcal{A}(d^{(\mu)},d^{(0)}) = a_{2} \cdot \frac{2}{(\mu+5)(\mu+4)(\mu+3)(\mu+2)}
\end{equation}
and if $\nu \ge 1$
\begin{equation}
  \mathcal{A}(d^{(\mu)},d^{(\nu)})
  =  a_{2} \cdot  \frac{\mu! \nu!}{(\mu+\nu+5)!}
      \left(2 \binom{\mu+\nu+2}{\nu+1} +
      \binom{\mu+\nu+2}{\nu} -
     \nu - 3 \right) \ .
  \label{eq:deuxieme}
\end{equation}
\end{auv}
\pr By Lemma 8
\begin{equation}
\begin{split}
  T_{\mu,\nu;p}(t) & = a_{2} \cdot \left(  2C(\mu,\nu) l^{\mu+\nu+4}
  + \sum_{k=0}^{\nu-1} \binom{\nu}{k} (\log^{\nu-k} p) l^{\mu+k+4}
    C(\mu,k) \right) \\
  & + O \left( \frac{l^{\mu+\nu+4}}{p} + l^{\mu+\nu+3} \right)
\end{split}
\end{equation}
where $l = \log t$ and $C(\mu,k)$ is defined by~(\ref{eq:cuk}).
Hence by the definition~(\ref{eq:Aab}) $\mathcal{A}(d^{(\mu)},d^{(0)})$ equals
\[
    a_{2} \cdot 2 \frac{\mu!}{(\mu+4)!}
     \left( \binom{\mu+2}{\mu+1} - 1 \right) \frac{(\mu+4)!}{(\mu+5)!}
   = a_{2} \cdot \frac{2}{(\mu+5)(\mu+4)(\mu+3)(\mu+2)}
\]
and if $\nu \ge 1$, $\mathcal{A}(d^{(\mu)},d^{(\nu)})$ equals
\begin{equation}
    a_{2} \cdot \left( 2 \frac{C(\mu,\nu)}{\mu+\nu+5}
    + \sum_{k=0}^{\nu-1} \left( \frac{\nu!}{k!(\nu-k)!} \right)
      \left( \frac{(\nu-k)!(\mu+k+4)!}{(\mu+\nu+5)!} \right)
      C(\mu,k) \right) \ .
      \label{eq:dieu}
\end{equation}
By~(\ref{eq:lampe}) the sum in~(\ref{eq:dieu}) is
\begin{equation}
 \frac{\mu! \nu!}{(\mu+\nu+5)!}
      \sum_{k=0}^{\nu-1}  \left(  \binom{\mu+k+2}{\mu+1}  -1
  \right)
 =  \frac{\mu! \nu!}{(\mu+\nu+5)!}
      \left( \binom{\mu+\nu+2}{\nu} - \nu - 1 \right) \ .
 \label{eq:maudite}
\end{equation}
Therefore~(\ref{eq:dieu}) and~(\ref{eq:maudite})
imply~(\ref{eq:deuxieme}).  \\

We summarize with a table of values of $\mathcal{A}(d^{(\mu)},d^{(\nu)})$.
In Table 2 the first column is the pair of sequences
$(d^{(\mu)},d^{(\nu)})$, the second column is the main term
of $T_{\mu,\nu;p}(t)$ as in~(\ref{eq:tuvp}),
and the third column is $\mathcal{A}(d^{(\mu)},d^{(\nu)})$
as computed by~(\ref{eq:deuxieme}).
Here we use the notation $l = \log t$ and $u = \log p$.
\begin{center}
{\bf Table 2}
\end{center}
\begin{center}
\begin{tabular}{|c|c|c|}
\hline
 $(d^{(\mu)},d^{(\nu)})$ & main term of $T_{\mu,\nu;p}(t)$
 & $\mathcal{A}(d^{(\mu)},d^{(\nu)})$ \\
\hline
\hline
$ (d,d) $ & $a_{2} \cdot \frac{1}{12} l^{4}$
            & $a_{2} \cdot \frac{1}{60}$  \\
\hline
$ (d^{(1)},d) $ & $ a_{2} \cdot \frac{1}{30} l^{5}$
                & $a_{2} \cdot  \frac{1}{180}$  \\
\hline
$ (d^{(2)},d) $ & $ a_{2} \cdot \frac{1}{60} l^{6}$
                  & $a_{2} \cdot  \frac{1}{420}$   \\
\hline
$ (d,d^{(1)}) $ & $a_{2} \cdot (\frac{1}{30}l^{5}+ \frac{1}{24}l^{4} u )
             $ &  $a_{2} \cdot \frac{1}{144}$  \\
\hline
$ (d^{(1)},d^{(1)}) $ & $ a_{2} \cdot (\frac{1}{72}l^6 + \frac{1}{60}l^5
 u) $ & $a_{2} \cdot \frac{1}{420} $ \\
\hline
$ (d^{(2)},d^{(1)}) $ & $ a_{2} \cdot  (\frac{1}{140}l^7 + \frac{1}{120}l^6
 u)$ & $a_{2} \cdot \frac{1}{960}$ \\
\hline
$ (d,d^{(2)}) $ & $ a_{2} \cdot  ( \frac{1}{60}l^6 + \frac{1}{30}l^5 u
                      + \frac{1}{24}l^4 u^2) $
                  & $a_{2} \cdot \frac{1}{280}$ \\
\hline
$ (d^{(1)},d^{(2)}) $ & $ a_{2} \cdot (\frac{1}{140}l^7 + \frac{1}{72} l^6 u
                     + \frac{1}{60} l^5 u^2) $
                     & $a_{2} \cdot \frac{5}{4032}$ \\
\hline
$ (d^{(2)},d^{(2)}) $ & $ a_{2} \cdot (\frac{19}{5040}l^8 + \frac{1}{140}l^7 u
                      + \frac{1}{120}l^6 u^2)$
                      & $a_{2} \cdot \frac{5}{9072}$ \\
\hline
\end{tabular}
\end{center}

Using the previous table we can compute $\mathcal{A}(a,b)$ for the
remainder of the sequences we require for our calculation.
\newtheorem{abc}[gonek]{Lemma}
\begin{abc}
We have
\begin{equation}
  \mathcal{A}(\e,d) = a_{2} \cdot \frac{1}{630} \ , \
  \mathcal{A}(d,\e) = a_{2} \cdot \frac{1}{420} \ , \
  \mathcal{A}(\e,d^{(1)}) = a_{2} \cdot \frac{1}{1440} \ , \
\end{equation}
\begin{equation}
  \mathcal{A}(d^{(1)},\e) = a_{2} \cdot \frac{17}{20160}  \ , \ and
  \ \mathcal{A}(\e,\e) = a_{2} \cdot \frac{23}{90720} \ .
\end{equation}
\end{abc}
\pr We employ the following notation:
if $A(t) = \sum_{n \le t} a_{n}$ and $j \in \mathbb{Z}_{\ge 0}$
then we define the operator $\mathcal{L}^{j}$ by
$(\mathcal{L}^{j} A)(t) = \sum_{n \le t} (\log^{j} n ) \ a_{n}$.
Note that if $A(t) = \sum_{n \le t} a_{n} = \alpha \log^{N} t
+ O( \log^{N-1} t )$
then partial summation implies
\begin{equation}
   (\mathcal{L}^{j} A) (t) = \alpha  \frac{N}{N+j}
   \log^{N+j} t + O( \log^{N+j-1} t ) \ .
   \label{eq:psid}
\end{equation}
By~(\ref{eq:froid}) we have the
identities
\begin{equation}
   \e(n) = \log n \ d^{(1)}(n) - d^{(2)}(n) \ , \
   \e(pn) = (\log p + \log n) d^{(1)}(pn) - d^{(2)}(pn) \ .
   \label{eq:enep}
\end{equation}
We begin with one example.  It follows from~(\ref{eq:enep})
that
\[
  d(n) \e(pn)
  = (\log p) \, d(n) d^{(1)}(pn) + (\log n) \,
  d(n) d^{(1)}(pn) - d(n)d^{(2)}(pn)
\]
and hence
\begin{equation}
   T_{d,\alpha;p}(t) = (\log p) \, T_{0,1;p}(t) +
   (\mathcal{L}T_{0,1;p})(t) - T_{0,2;p}(t) \ .
   \label{eq:tdap}
\end{equation}
By Table 2,~(\ref{eq:psid}), and~(\ref{eq:tdap}) we
derive
\begin{equation}
  T_{d,\e;p}(t) = u T_{0,1;p}(t) +
  (\mathcal{L}T_{0,1;p})(t) - T_{0,2;p}(t)
  = a_{2} \cdot \left( \frac{l^{6}}{90} + \frac{l^{5}u}{30} \right)
     + O \left( l^{5} + \frac{l^{6}}{p} \right)
\end{equation}
where $l = \log t$ and $u = \log p$.
In a similar fashion we compute
\begin{equation}
\begin{split}
    & T_{\alpha,d;p}(t) = (\mathcal{L}T_{1,0;p})(t) - T_{2,0;p}(t) \ ,
    \\
    & T_{\alpha,d^{(1)};p}(t) = (\mathcal{L}T_{1,1;p})(t) - T_{2,1;p}(t) \ ,
    \\
    & T_{d^{(1)},\alpha;p}(t) = (\log p) \, T_{1,1;p}(t)
                      + (\mathcal{L}T_{1,1;p})(t) - T_{1,2;p}(t) \ , \ and  \\
    & T_{\alpha,\alpha;p}(t) = (\log p) \, (\mathcal{L}T_{1,1;p})(t)
                 + (\mathcal{L}^{2}T_{1,1;p})(t)
                 - (\mathcal{L}T_{1,2;p})(t) \\
                 & - (\log p) \, T_{2,1;p}(t)
                 -   (\mathcal{L}T_{2,1;p})(t)
                  + T_{2,2;p}(t) \ . \\
    \label{eq:iden}
\end{split}
\end{equation}
Thus Table 2,~(\ref{eq:psid}), and~(\ref{eq:iden}) imply
\begin{equation}
   T_{a,b;p}(t) =
   \sum_{i+j=A} c_{ij}  (\log^{i} t )
   ( \log^{j} p) + O \left( \frac{\log^{A} t}{p} + \log^{A-1} t \right) \
\end{equation}
for the aforementioned sequences ($a$,$b$) and appropriate constants
$c_{ij},A$.
In summary, we obtain
\begin{center}
{\bf Table 3}
\end{center}
\begin{center}
\begin{tabular}{|c|c|c|}
\hline
 $(a,b)$ & main term of $T_{a,b;p}(t)$
& $\mathcal{A}(a,b)$ \\
\hline
\hline
$ (\alpha,d) $ & $  a_{2} \cdot \frac{1}{90} l^6 $
  & $ a_{2} \cdot \frac{1}{630}$ \\
\hline
$ (d,\alpha) $ & $
 a_{2} \cdot (\frac{1}{90}l^6 + \frac{1}{30}l^5 u ) $
  & $ a_{2} \cdot \frac{1}{420}$ \\
\hline
$ (\alpha,d^{(1)}) $ & $
              a_{2} \cdot  (\frac{1}{210}l^7 + \frac{1}{180}l^6 u)$
  & $ a_{2} \cdot \frac{1}{1440}$ \\
\hline
$ (d^{(1)},\alpha) $ & $
              a_{2} \cdot ( \frac{1}{210}l^7 + \frac{1}{72}l^6 u) $
  & $ a_{2} \cdot \frac{17}{20160}$ \\
\hline
$ (\alpha,\alpha) $
 & $ a_{2} \cdot \left( \frac{17}{10080} l^8 + \frac{1}{210} l^7u
 \right)$
 & $ a_{2} \cdot \frac{23}{90720}$ \\
\hline
\end{tabular}
\end{center}
\subsection{Proof of Corollary 1}
\pr We first evaluate $S_{\alpha} = I(\alpha,\alpha;T)$. We have
$\alpha(n) \ll d(n) \log^{2}(n)$,
\[
  \sum_{n \le t} \frac{\alpha(n)^{2}}{n} = a_{2}
   \cdot \frac{17}{20160} \log^{8} t +
    O( \log^{7} t) \ , \ and
\]
\[
   \sum_{n \le t} \frac{\alpha(n)\alpha(pn)}{n} =
    a_{2} \cdot
   \left( \frac{17}{10080} \log^8 t + \frac{1}{210} \log^7 t
    \log p  \right)
   + O \left( \log^{7} t
              + \frac{\log^{8} t}{p} \right)
\]
by Tables 1 and 3. Moreover,
$\mathcal{A}(\e,\e) = a_{2} \cdot  \frac{23}{90720}$
by Lemma 10. Applying Lemma 5, we deduce
\begin{equation}
 S_{\e} =  \frac{T L^{9}}{2 \pi}
          \left( a_{2} \cdot \frac{17}{20160} - 2\mathcal{A}(\e,\e)
          \right)
       \sim \frac{61}{60480 \pi^{3}} T L^{9}
       = \frac{a_{2}TL^{9}}{2 \pi} \frac{61}{181440}
  \label{eq:ciel}
\end{equation}
with an error term $O(TL^{8} \log L)$.
Next we consider
\begin{equation}
  S_{\beta} = \sum_{0 < \gamma < T} D_{\beta_{\gamma}}(\rho)
                                     D_{\beta_{\gamma}}(1-\rho)
            = \sum_{0 < \gamma < T} \sum_{1 \le m,n \le
                                     \frac{\gamma}{2\pi}}
                                    \frac{\beta_{\gamma}(m)
                                     \beta_{\gamma}(n)}
                                    {m^{\rho}n^{1-\rho}} \ .
  \label{eq:ame}
\end{equation}
Before evaluating $S_{\beta}$, we require some notation.
For $N \in \mathbb{Z}_{\ge 0}$, $a(n)$ and $b(n)$ sequences,
define
\[
       I_{N}(a,b;T) =
       \sum_{0 < \gamma < T} \log^{N}
      \left( \frac{\gamma}{2 \pi} \right)
       D_{a}(\rho)D_{b}(1-\rho)  \ .
       \label{eq:Skab}
\]
Notice that $I_{0}(a,b;T) = I(a,b;T)$ of~(\ref{eq:iabt}). Observe
that if $a(n)$ and $b(n)$ are real sequences then $I_{N}(b,a;T) =
\overline{I_{N}(a,b;T)}$ by a consideration similar
to~(\ref{eq:vache}). Applying the identity $\beta_{t}(m) = l^{2}
d(m) - 2l d^{(1)}(m) + \e(m)$ where $l = \log( \frac{t}{2 \pi})$
we obtain
\begin{equation}
\begin{split}
 \f(m) \f(n)
 & = l^{4} d(m) d(n) - 2l^{3}d(m) d^{(1)}(n) + l^{2}d(m)\e(n) \\
 & - 2l^{3} d^{(1)}(m) d(n) + 4l^{2} d^{(1)}(m) d^{(1)}(n)
   - 2l d^{(1)}(m)\e(n) \\
 & + l^{2} \e(m) d(n) - 2l \e(m) d^{(1)}(n) + \e(m)\e(n) \ . \\
 \label{eq:souris}
\end{split}
\end{equation}
Inserting~(\ref{eq:souris}) in~(\ref{eq:ame}) we obtain
\begin{equation}
\begin{split}
  S_{\beta}
          & = I_{4}(d,d;T) + 4I_{2}(d^{(1)},d^{(1)};T)
                 + I(\e, \e; T)   \\
          &     - 4 \mathrm{Re}(I_{3}(d,d^{(1)};T))
               + 2 \mathrm{Re}( I_{2}(d,\e;T))
               -4 \mathrm{Re}( I_{1}(d^{(1)},\e;T) ) \ .  \\
  \label{eq:enfer}
\end{split}
\end{equation}
Note that if $I(a,b;T) = c_{1}T \log^{M} T + O(T \log^{M-1} T)$
then partial summation implies $ I_{N}(a,b;T) = I(a,b;T)(\log^{M}
T + O(\log^{M-1} T))$. In an analogous calculation to that of
$S_{\alpha}$, we derive by Lemma 5, Tables 1-3, and the partial
summation identity the following:
\begin{equation}
  I(d, d;T) \sim \frac{a_{2}T L^{9}}{2 \pi}
        \left( \frac{1}{24} - 2 \cdot \frac{1}{60}
          \right)
        = \frac{a_{2}T L^{9}}{2 \pi} \cdot \frac{1}{120} \ ,
  \label{eq:mouche}
\end{equation}
\begin{equation}
  I(d^{(1)}, d^{(1)};T)  \sim \frac{a_{2}T L^{9}}{2 \pi}
        \left( \frac{1}{144} - 2 \cdot \frac{1}{420}
          \right)
        = \frac{a_{2}T L^{9}}{2 \pi} \cdot \frac{11}{5040} \ ,
\end{equation}
\begin{equation}
  I_{3}(d,d^{(1)};T)  \sim \frac{a_{2} T L^{9}}{2 \pi}
         \left( \frac{1}{60} - \frac{1}{144} - \frac{1}{180}
         \right)
        = \frac{a_{2} T L^{9}}{2 \pi} \cdot \frac{1}{240} \ ,
\end{equation}
\begin{equation}
  I_{2}(d,\e;T)  \sim \frac{a_{2} T L^{9}}{2 \pi}
        \left( \frac{1}{180} - \frac{1}{420} - \frac{1}{630}
        \right)
        = \frac{a_{2} T L^{9}}{2 \pi} \cdot \frac{1}{630} \ ,
\end{equation}
\begin{equation}
  I_{1}(d^{(1)},\e;T)  \sim \frac{a_{2} T L^{9}}{2 \pi}
          \left( \frac{1}{420} - \frac{17}{20160}
          - \frac{1}{1440} \right)
        = \frac{a_{2} T L^{9}}{2 \pi} \cdot \frac{17}{20160}
  \label{eq:gueppe}
\end{equation}
where each of these holds with an error term $O(TL^{8} \log L)$.
\noindent By~(\ref{eq:enfer}),~(\ref{eq:ciel})
and~(\ref{eq:mouche}) -~(\ref{eq:gueppe}) we have
\[
  S_{\beta} \sim \frac{a_{2}T L^{9}}{2 \pi}
       \left( \frac{1}{120}
       +  4 \cdot \frac{11}{5040} + \frac{61}{181440}
       - 4 \cdot \frac{1}{240}
       + 2 \cdot \frac{1}{630} - 4 \cdot \frac{17}{20160} \right) \ .
\]
This simplifies to $S_{\beta} =  \frac{97}{60480 \pi^3} T L^{9} +
O(T L^{8} \log L)$.
\subsection{Proof of Theorem 2}
\pr By Tables 1 and 2 we have
\begin{equation}
  T_{d,d}(t)
  = a_{2} \cdot \frac{1}{24} \log^{4} t + O(\log^{3} t) \ , \
  T_{d,d;p}(t)
  = a_{2} \cdot \frac{1}{12} \log^{4} t + O(\log^{3} t) \ .
  \label{eq:form}
\end{equation}
Therefore by Lemma 3
\begin{equation}
   I(d,d;T,\delta) =
   \frac{T}{2 \pi} \left( \frac{a_{2}}{24} L^{5}
                          - 2 \mathrm{Re}
                   \left( M \left(d,d;\frac{T}{2 \pi},
   \delta \right) \right) \right) + O(TL^{4} \log L)
   \label{eq:iddt}
\end{equation}
where $\delta = \frac{\lambda}{L}$. Since we have~(\ref{eq:form})
an application of Lemma 4 yields
\begin{equation}
   M \left(d,d;\frac{T}{2 \pi},
   \delta \right)
   = L^{5} \sum_{k=0}^{\infty} \frac{(i \delta L)^{k}}{k!}
   \frac{a_{2}}{12} \frac{k! 4!}{(5+k)!} + O(L^{4}) \ .
   \label{eq:mddt}
\end{equation}
Thus~(\ref{eq:iddt}) and~(\ref{eq:mddt}) imply
\begin{equation}
    I(d,d;T,\delta)
    = \frac{a_{2}TL^{5}}{2 \pi}
    \left( \frac{1}{5!} - 4 \sum_{j \ge 1} \frac{(-1)^{j}
    \lambda^{2j}}{(5+2j)!} \right)
     + O(TL^{4} \log L)
\end{equation}
and we are finished.

\noindent D\'epartement de Math\'ematiques et de statistique,
Universit\'e de Montr\'eal,
CP 6128 succ Centre-Ville,
Montr\'eal, QC, Canada  H3C 3J7 \\
EMAIL: nathanng@dms.umontreal.ca


\begin{thebibliography}{99}
\bibitem{C}
J.B. Conrey, {\em The fourth moment of derivatives of the Riemann
zeta-function}, Quart. J. Math. Oxford (2), \textbf{39} (1988),
21-36.

\bibitem{CGG}
J.B. Conrey, A. Ghosh, S.M. Gonek, {\em Simple zeros of the Riemann
  zeta function}, Proc. London Math Soc. (3) \textbf{76} (1998),
no. 3, 497-522.

\bibitem{G1}
S.M. Gonek, {\em Mean values of the Riemann zeta-function and its derivatives},
Invent. Math. \textbf{75} (1984), 123-141.

\bibitem{G2}
S.M. Gonek, {\em An explicit formula of Landau and its applications
to the theory of the zeta function}, Contemp. Math \textbf{143}, 1993,
395-413.

\bibitem{G3}
S.M. Gonek, {\em On negative moments of the Riemann zeta-function},
Mathematika \textbf{36} (1989), 71-88.

\bibitem{GKP}
R.L. Graham, D.E. Knuth, O. Patashnik, {\em Concrete Mathematics,
second edition}, Addison-Wesley, New York, 1994.

\bibitem{Ha}
R.R. Hall, {\em The behaviour of the Riemann zeta-function on the
critical line}, Mathematika \textbf{46} (1999), 281-313.

\bibitem{He}
D. Hejhal, {\em On the distribution of log$\zeta^{'}(1/2+it)$}, Number
Theory, Trace Formulas, and Discrete Groups, K.E. Aubert, E. Bombieri
and D.M. Goldfeld, eds., Proceedings of the 1987 Selberg Symposium,
(Academic Press, 1989), 343-370.

\bibitem{HKO}
C.P. Hughes, J.P. Keating, and Neil O'Connell,
{\em Random matrix theory and the derivative of the Riemann zeta
function}, Proceedings of the Royal Society: A \textbf{456} (2000),
2611-2627.

\bibitem{In}
A.E. Ingham, {\em Mean-value theorems in the theory of the Riemann
zeta function}, P.L.M.S. (2), \textbf{27} (1926), 273-300.

\bibitem{KS}
J.P. Keating and N.C. Snaith, {\em Random matrix theory and
$\zeta(1/2 + it)$}, Communications in Mathematical Physics
\textbf{214} (2000), 57-89.

\bibitem{Mu}
J. Mueller, {\em On the difference between consecutive zeros of the
Riemann zeta function}, J. Number Theory \textbf{14} (1982), 327-331.

\bibitem{Ng}
N. Ng, {\em Limiting distributions and zeros of Artin $L$-functions},
Ph. D. Thesis, University of British Columbia, fall 2000.

\bibitem{Sh}
P. Shiu,  {\em A Brun-Titchmarsh theorem for multiplicative
functions}, J. Reine Angew. Math. \textbf{313} (1980), 161-170.

\bibitem{Te}
G. Tenenbaum, {\em Introduction to analytic and probabilistic number
  theory}, Cambridge University Press, Cambridge, 1995.

\bibitem{Ti}
E.C. Titchmarsh, {\em The theory of the Riemann zeta function, second
edition}, Oxford University Press, New York, 1986.
\end{thebibliography}
\end{document}